\documentclass[10pt]{amsart}
\usepackage{mathptmx}
\usepackage{amsmath}
\usepackage{amssymb}
\usepackage{array}
\usepackage{geometry}
\usepackage[bookmarks=true,colorlinks=true, pdfstartview=FitV, linkcolor=black, citecolor=blue, urlcolor=black]{hyperref}

\usepackage{color}
\definecolor{DarkRed}{rgb}{0.55,.00,0.2}
\definecolor{DarkGrey}{rgb}{0.35,.35,0.35}

\theoremstyle{definition}

\theoremstyle{remark}

\numberwithin{equation}{section}

\begin{document}

\title{  New index transforms  with the product of Bessel functions}

\author{S. Yakubovich}
\address{Department of Mathematics, Faculty of Sciences,  University of Porto,  Campo Alegre str.,  687; 4169-007 Porto,  Portugal}
\email{ syakubov@fc.up.pt}


\date{\today}
\maketitle

\markboth{\rm \centerline{ S.  Yakubovich}}{}
\markright{\rm \centerline{Index Transforms  With the Product of Bessel Functions }}

\begin{abstract}  New index transforms are  investigated, which contain as the kernel products of the Bessel and modified Bessel functions.   Mapping properties and  invertibility in  Lebesgue  spaces are studied  for these operators.  Relationships  with the Kontorovich-Lebedev and Fourier cosine transforms are  established.   Inversion theorems are proved.   As an application, a solution of the initial value problem for the  fourth order partial differential equation, involving the Laplacian  is presented. 
\end{abstract}

\section{Introduction and preliminary results}

The main objects of the present paper are  the following operators of the index transforms \cite{yak}
$$(Ff) (\tau) = {1\over \sinh(\pi\tau /2) }\int_0^\infty  K_{i\tau}(2\sqrt {2x}) \ {\rm Im} \left[ J_{i\tau} (2\sqrt {2x})\right] f(x)dx,   \quad    \tau \in \mathbb{R} \backslash \{0\}, \eqno(1.1)$$
$$ (Gg) (x) = \int_{-\infty}^\infty    K_{i\tau}(2\sqrt {2x}) \ {\rm Im} \left[ J_{i\tau} (2\sqrt {2x}) \right] \frac{ g(\tau)} {\sinh(\pi\tau /2)} \  d\tau,   \quad    x \in \mathbb{R}_+. \eqno(1.2)$$
Here  $i$ is the imaginary unit and ${\rm Im}$ denotes the imaginary part of a complex-valued function.  The convergence of the integrals (1.1), (1.2) will be clarified below. Functions  $J_\mu(z), K_\mu(z),\  \mu \in \mathbb{C}$ \cite{erd}, Vol. II  are,  correspondingly,  the Bessel and  modified Bessel or Macdonald functions,   which satisfy  the differential equations
$$  z^2{d^2u\over dz^2}  + z{du\over dz} +  (z^2- \mu^2)u = 0,\eqno(1.3)$$
$$  z^2{d^2u\over dz^2}  + z{du\over dz} -  (z^2+ \mu^2)u = 0,\eqno(1.4)$$
respectively.   It has the asymptotic behaviour 
$$J_\mu(z)= \sqrt{2\over \pi z} \cos\left(z- {\pi\over 4} (2\mu +1)\right) [1+
O(1/z)], \quad   K_\mu(z) = \left( \frac{\pi}{2z} \right)^{1/2} e^{-z} [1+
O(1/z)], \qquad z \to \infty,\eqno(1.5)$$
and near the origin
$$J_\mu(z)= O(z^\mu), \quad z^\mu K_\mu(z) = 2^{\mu-1} \Gamma(\mu) + o(1), \ z \to 0,\eqno(1.6)$$
$$K_0(z) = -\log z + O(1), \ z \to 0. \eqno(1.7)$$
The Macdonald function can be represented by the integral 
$$K_{\mu}(z)= \int_0^\infty e^{- z\cosh u}\cosh (\mu u)du,\   {\rm Re}\  z >0,  \   \mu \in \mathbb{C}.\eqno(1.8)$$
Concerning the kernel of the transform (1.1) our key representation will be in terms of the Mellin-Barnes integral, giving by relation (8.4.23.11) in \cite{prud}, Vol. III, namely,
$$  { K_{i\tau}(2\sqrt {2x})  \over \sinh(\pi\tau /2) }\ {\rm Im} \left[ J_{i\tau} (2\sqrt {2x} )\right] = - {1\over 16\pi i \sqrt\pi} \int_{\gamma-i\infty}^{\gamma+i\infty} \frac{ \Gamma(s/2)\Gamma((s+i\tau)/2) \Gamma((s-i\tau)/2)} {\Gamma((1-s)/2)} x^{-s} ds,\eqno(1.9)$$
where $\  x >0,\  \tau \in \mathbb{R}\backslash \{0\} ,\  \gamma >0,\  \Gamma(z)$ is Euler's gamma -function \cite{erd}, Vol. I.    The theory of the index transforms can be found in the author monograph \cite{yak} (see also \cite{yal}).   The familiar example is the Kontorovich-Lebedev transform
$$(KL f)(\tau)= \int_0^\infty K_{i\tau}(x) f(x) dx.\eqno(1.10)$$
Our method of investigation is based on the theory of the Mellin transform in the Lebesgue spaces \cite{tit}.   In fact, we define the Mellin transform  as 
$$f^*(s)= \int_0^\infty f(x) x^{s-1} dx\eqno(1.11)$$
and its inverse, accordingly,
 $$f(x)= {1\over 2\pi i}  \int_{\nu- i\infty}^{\nu+i\infty} f^*(s)  x^{-s} ds.\eqno(1.12)$$
 Integrals (1.11), (1.12) are convergent, for instance, in mean in the weighted $L_p$-spaces,  $1<p\le 2$ and 
 the Parseval equality takes place \cite{tit}
$$\int_0^\infty f(x) g(x) dx= {1\over 2\pi i} \int_{\nu- i\infty}^{\nu+i\infty} f^*(s) g^*(1-s) ds.\eqno(1.13)$$
It is important for us for further investigation to obtain an integral representation of the kernel in (1.1) in terms of the Fourier cosine transform \cite{tit}
$$(F_cf)(x) = \sqrt{2\over \pi} \int_0^\infty f(t) \cos(xt) dt.\eqno(1.14)$$
Precisely, we prove the following 

{\bf Lemma 1}. {\it Let $ x >0, \tau \in \mathbb{R} \backslash \{0\}$. Then
$${K_{i\tau}(2\sqrt {2x})  \over \sinh(\pi\tau /2) }  {\rm Im} \left[ J_{i\tau} (2\sqrt {2x})\right]   =  -  {2\over \pi }  \int_0^\infty   \cos(\tau u) \ {\rm Re}  \left[ K_0\left( 4 e^{\pi i/4} \sqrt{ x\cosh u} \right) \right]\ du,\eqno(1.15)$$
where ${\rm Re}$ denotes the real part of a complex-valued function. }

\begin{proof} Indeed, taking the Fourier cosine transform (1.14) from both sides of the equality (1.9), we change the order of integration owing to Fubini's theorem.  In fact,  employing the reciprocal equalities (see formula (1.104) in \cite{yak})  via the Fourier cosine transform (1.14) 
$$\int_0^\infty  \Gamma\left({s\over 2} + {i\tau\over 2}\right)  \Gamma\left({s\over 2} - {i\tau\over 2}\right)  \cos(u \tau) d\tau
= {\pi\over 2^{s-1}}  {\Gamma(s) \over \cosh^s u},\ {\rm Re}\ s > 0,\eqno(1.16)$$
$$  \Gamma\left({s\over 2} + {i\tau\over 2}\right)  \Gamma\left({s\over 2} - {i\tau\over 2}\right)  
=   { \Gamma(s)  \over 2^{s-2}}  \int_0^\infty   {\cos(\tau u)  \over \cosh^s u} \ du,\eqno(1.17)$$ 
one can integrate twice by parts in the integral (1.17), showing the uniform estimate
$$\left| \Gamma\left({s\over 2} + {i\tau\over 2}\right)  \Gamma\left({s\over 2} - {i\tau\over 2}\right)  \right|
\le {|\Gamma(s+1)|\over \tau^2} \left[ c_1+ c_2 |s|\  \right],\  {\rm Re}\ s > 0,\  \tau \in  \mathbb{R} \backslash \{0\} ,\eqno(1.18)$$
where $c_1, c_2$ are absolute positive constants.  Hence with the use of  the Stirling asymptotic formula for the gamma-function \cite{erd}, Vol. I  and the elementary inequality 
$$ \left| \Gamma\left({s\over 2} + {i\tau\over 2}\right)  \Gamma\left({s\over 2} - {i\tau\over 2}\right)  \right| \le |\Gamma (s)|
B(\gamma/2,\gamma/2)\eqno(1.19)$$
where  $B(a,b)$ is Euler's beta-function \cite{erd}, Vol. I,  it guarantees the absolute convergence of the iterated integral 
$$\int_0^\infty   \int_{\gamma-i\infty}^{\gamma+i\infty} \left|\frac{ \Gamma(s/2)\Gamma((s+i\tau)/2) \Gamma((s-i\tau)/2)} {\Gamma((1-s)/2)}  ds\right| d\tau   < \infty.$$
Therefore, (1.9) and (1.16) yield 
$$ \int_0^\infty   {\cos(\tau u) \over \sinh(\pi\tau /2) } K_{i\tau}(2\sqrt {2x})\  {\rm Im} \left[ J_{i\tau} (2\sqrt {2x})\right]  d\tau = - {1\over 8i \sqrt\pi} \int_{\gamma-i\infty}^{\gamma+i\infty} \frac{ \Gamma(s) \Gamma(s/2)} {\Gamma((1-s)/2)} (2x\cosh u)^{-s} ds.\eqno(1.20)$$
Moreover, appealing to the duplication formula for the gamma-function \cite{erd}, Vol. I and making a simple substitution, the latter Mellin-Barnes integral  can be written as 
$$ - {1\over 8i \sqrt\pi} \int_{\gamma-i\infty}^{\gamma+i\infty} \frac{ \Gamma(s) \Gamma(s/2)} {\Gamma((1-s)/2)} (2x\cosh u)^{-s} ds = - {1\over 8 \pi i} \int_{\gamma-i\infty}^{\gamma+i\infty} \frac{ \Gamma^2(s) \Gamma(s+ 1/2)} {\Gamma(1/2-s)} (x^2\cosh^2 u)^{-s} ds.\eqno(1.21)$$
Meanwhile, considering $0 < \gamma < 1/2$, we have the value of the relatively convergent integral (cf. \cite{tit}, Section 7.3)
$$ \frac{ \Gamma(s) } {\Gamma(1/2-s)} = {2^{1-2s}\over \sqrt \pi} \int_0^\infty \cos t\   t^{2s-1} dt
= -  {2^{1-2s} (2s-1) \over \sqrt \pi} \int_0^\infty \sin  t\   t^{2(s-1)} dt,  $$
where the latter integral is obtained after the integration by parts and, evidently,  converges absolutely for $0 < \gamma < 1/2$.   Hence, substituting it in the right-hand side of (1.21) and changing the order of integration, we find with the aid of relation (8.4.23.1) in \cite{prud}, Vol. III 
$$- {1\over 8\pi i} \int_{\gamma-i\infty}^{\gamma+i\infty} \frac{ \Gamma^2(s) \Gamma(s+ 1/2)} {\Gamma(1/2-s)} (x^2\cosh^2 u)^{-s} ds$$
$$ =  {1\over 4\pi i \sqrt \pi} \int_0^\infty {\sin  t\over t^2} \  \int_{\gamma-i\infty}^{\gamma+i\infty}  \left[2 \Gamma(1+s) \Gamma(s+ 1/2) -  \Gamma(s) \Gamma(s+ 1/2) \right] \left({2x\cosh u\over t}\right)^{-2s} ds dt $$
$$= {1 \sqrt { 2x\cosh u} \over  \sqrt \pi} \int_0^\infty \  K_{1/2} \left( {4x\cosh u\over t}\right)
\left[   {4x\cosh u\over t}  -  1 \right] \   {\sin  t\over t^{5/2} } dt.$$
Since
$$K_{1/2}(z)= \sqrt{{ \pi\over 2z}}  e^{-z},$$
we get with the simple substitution 
$$- {1\over 8 \pi i} \int_{\gamma-i\infty}^{\gamma+i\infty} \frac{ \Gamma^2(s) \Gamma(s+ 1/2)} {\Gamma(1/2-s)} (x^2\cosh^2 u)^{-s} ds =  {1\over 2}   \int_0^\infty \  e^{- 4x \cosh( u) / t} \left[   {4x  \cosh u\over t}    -  1 \right] \  { \sin t \over t^2} \   dt.$$
We calculate the latter integral, appealing to relation (2.5.37.1) in \cite{prud}, Vol. I. In fact, taking in mind the identity for  Macdonald's  functions  \cite{erd}, Vol. II
$$K_{\mu+1}(z)- K_{\mu-1}(z)= {2\mu\over z} K_\mu(z),$$
we deduce
$$ {1\over 2}  \int_0^\infty \  e^{- 4x \cosh( u) / t} \left[   {4x  \cosh u\over t}    -  1 \right] \  { \sin t \over t^2} \   dt 
= - {1\over 2}  \left[ K_2\left( 4 e^{\pi i/4} \sqrt{ x\cosh u} \right) + K_2\left( 4 e^{- \pi i/4} \sqrt{ x\cosh u} \right) \right]$$
$$+ {1\over 4i \sqrt{x\cosh u}} \left[   e^{\pi i/4}  K_1\left( 4 e^{\pi i/4} \sqrt{ x\cosh u} \right) -  e^{- \pi i/4}  K_1 \left( 4 e^{- \pi i/4} \sqrt{ x\cosh u} \right) \right] $$ 
$$= - {\rm Re}  \left[ K_0\left( 4 e^{\pi i/4} \sqrt{ x\cosh u} \right) \right].$$ 
Thus combining with (1.20), (1.21), we find the value of the index integral 
$$  \int_0^\infty   {\cos(\tau u) \over \sinh(\pi\tau /2) } K_{i\tau}(2\sqrt {2x}) \ {\rm Im} \left[ J_{i\tau} (2\sqrt {2x})  \right] d\tau =  -   {\rm Re} \left[ K_0\left( 4 e^{\pi i/4} \sqrt{ x\cosh u} \right)\right],\  x >0, \ u \in \mathbb{R}.\eqno(1.22)$$ 
In the meantime, using (1.9),  (1.18), one can verify that for each $ x >0$ the kernel of (1.1) belongs to $L_1(\mathbb{R}_+) \cap L_2(\mathbb{R}_+)$.    Since (see (1.5), (1.7), (1.8) ) 
$$\left| K_0\left( 4 e^{\pm \pi i/4} \sqrt{ x\cosh u} \right)\right| = \left| \int_0^\infty  \exp \left( -  4 e^{\pm \pi i/4} \sqrt{ x\cosh u} \  \cosh t\right)  dt \right| $$
$$\le  \int_0^\infty  \exp \left( -  2  \sqrt{2 x\cosh u} \  \cosh t\right)  dt = K_0( 2  \sqrt{2 x\cosh u}), $$
the same is true for the right-hand side of (1.22). Hence taking the inverse Fourier cosine transform \cite{tit}, we come up with equality (1.15), completing the proof of Lemma 1. 
\end{proof}

{\bf Corollary 1}. {\it It has }
$$\lim_{\tau \to 0} \quad  {K_{i\tau}(2\sqrt {2x})  \over \sinh(\pi\tau /2) }  {\rm Im} \left[ J_{i\tau} (2\sqrt {2x})  \right] =  -   {2\over \pi}  \int_0^\infty   {\rm Re} \left[ K_0\left( 4 e^{\pi i/4} \sqrt{ x\cosh u} \right)  \right]  du,\  x >0.$$

{\bf Corollary 2}. { \it Let  $ x >0, \tau \in \mathbb{R} \backslash \{0\}$. Then the kernel in $(1.1)$  satisfies the following inequality 
$$\left|  {K_{i\tau}(2\sqrt {2x})  \over \sinh(\pi\tau /2) }  \ {\rm Im} \left[ J_{i\tau} (2\sqrt {2x}) \right] \right| 
 \le    {4\over \pi} \  e^{-\delta |\tau| } \  K^2_0\left( \cos \left( {\delta\over 2} \right) \sqrt{2 x \cos(\delta)} \right),\eqno(1.23)$$
where  $\delta \in [ 0, \pi/2)$.} 

\begin{proof}  Indeed, the integral  in the right-hand side of (1.15) can  be written as 
$$  - {1\over 2 \pi }  \int_{-\infty}^\infty   e^{i\tau u}   \left[ K_0\left( 4 e^{\pi i/4} \sqrt{ x\cosh u} \right)  +  K_0\left( 4 e^{- \pi i/4} \sqrt{ x\cosh u} \right)\right] du$$
and by the analytic property of the integrand and the absolute convergence of the integral one can move the contour along the open infinite horizontal strip $(i\delta - \infty,\ i\delta + \infty)$ with $\delta \in [ 0, \pi/2)$, i.e. 
$$ - {1\over 2 \pi}  \int_{-\infty}^\infty   e^{i\tau u} \    \left[ K_0\left( 4 e^{\pi i/4} \sqrt{ x\cosh u} \right)  +  K_0\left( 4 e^{- \pi i/4} \sqrt{ x\cosh u} \right) \right] du$$
$$ =  -  {1\over 2 \pi }  \int_{-\infty}^\infty   e^{i\tau (i\delta+  u)}  \  \left[ K_0\left( 4 e^{\pi i/4} \sqrt{ x\cosh (i\delta+ u)} \right)  + K_0\left( 4 e^{- \pi i/4} \sqrt{ x\cosh (i\delta + u)}\right) \right] du,\eqno(1.24)$$
where the main branch of the square root is chosen.   Hence denoting by 
$$z= \cosh (i\delta+ u) = \cos(\delta)\cosh u + i \sin(\delta)\sinh u= |z| e^{i \arg z }, $$
where $\arg z= \arctan \left( \tan(\delta) \tanh u\right) \in ( -\delta,\  \delta).$   Then 
$${\rm Re } \left[ 4\  e^{\ \pm \pi i/4} \sqrt{ x\cosh (i\delta+ u)} \right]  = 4 \sqrt{ x |z| } \cos \left( {\arg z\over 2} \pm {\pi\over 4}\right) \ge  4 \sqrt{ x \cos(\delta)\cosh u } \  \cos \left( {\arg z\over 2} \pm {\pi\over 4}\right) $$
 when   $\cos \left( {\arg z\over 2} \pm {\pi\over 4}\right) = {1\over \sqrt 2} \left[ \cos \left( \arg z /2 \right)  \mp   \sin  \left( \arg z /2 \right) \right] > 0$, i.e  $|\tan(\arg z /2)| < 1,$ which is true.  Moreover, we have 
 $$ \cos \left( {\arg z\over 2} \pm {\pi\over 4}\right)  >  \cos \left( {\delta\over 2}  + {\pi\over 4}\right).$$ 
 Therefore, returning to (1.15), (1.24), we obtain
 $$\left|  {K_{i\tau}(2\sqrt {2x})  \over \sinh(\pi\tau /2) } \ {\rm Im}  \left[ J_{i\tau} (2\sqrt {2x})  \right] \right| 
 \le   {2 e^{-\delta |\tau| } \over \pi }  \int_{0}^\infty   K_0\left( 4 \cos \left( {\delta\over 2}  + {\pi\over 4}\right) \sqrt{ x \cos(\delta)\cosh u }  \right)  du $$
 $$\le   {2 e^{-\delta |\tau| } \over \pi }  \int_{0}^\infty   K_0\left( 4 \cos \left( {\delta\over 2}  + {\pi\over 4}\right) \sqrt{ x \cos(\delta)} \cosh (u/2)\right)  du $$
 $$=   {4 e^{-\delta |\tau| } \over \pi }  \int_{0}^\infty  \int_{0}^\infty \exp\left(-  4 \cos \left( {\delta\over 2}  + {\pi\over 4}\right) \sqrt{ x \cos(\delta)} \cosh u \  \cosh t \right)  du\ dt $$
 $$\le  {4 e^{-\delta |\tau| } \over \pi }   K^2_0\left( \cos \left( {\delta\over 2} \right) \sqrt{ 2 x \cos(\delta)} \right),$$
 completing the proof of Corollary 2.   
 
 \end{proof}
 
 An immediate consequence of Lemma 1 is also 
 
 {\bf Corollary 3}. { \it Let  $ x >0$. Then for all $ \tau \in \mathbb{R} $ the inequality   
$$\left|  {K_{i\tau}(2\sqrt {2x})  \over \sinh(\pi\tau /2) }  \ {\rm Im} \left[ J_{i\tau} (2\sqrt {2x}) \right] \right| 
 \le    {4\over \pi}   \ K_0^2 \left( \sqrt {2x} \right) \eqno(1.25)$$
is fulfilled.} 
 
 \begin{proof}  We have
 $$ {2\over \pi } \left| \int_0^\infty   \cos(\tau u) \ {\rm Re}  \left[ K_0\left( 4 e^{\pi i/4} \sqrt{ x\cosh u} \right) \right]\ du\right|
 \le {2\over \pi} \int_0^\infty    K_0\left( 2 \sqrt{ 2 x\cosh u} \right) \ du $$
 $$\le {4\over \pi} \int_0^\infty    K_0\left( 2 \sqrt{ 2 x} \  \cosh u \right) \ du =  {4\over \pi} \int_0^\infty     \int_0^\infty    \exp \left( - \sqrt{ 2 x} \  \cosh u  \  \cosh v \right)$$
 $$\times  \   \exp \left( - \sqrt{ 2 x} \  \cosh u  \  \cosh v \right) \ du dv\le {4\over \pi}  \ K_0^2 \left( \sqrt {2x} \right) $$ 
  and the result follows. 

\end{proof} 

{\bf Remark 1}.  Inequality (1.25) is a particular case of the inequality (1.23) with $\delta=0$. 

Employing the Mellin-Barnes representation (1.9) of the kernel in (1.1), which we denote by 
$$\Psi_{\tau}(x)=  { K_{i\tau}(2\sqrt {2x})  \over \sinh(\pi\tau /2) }\ {\rm Im} \left[ J_{i\tau} (2\sqrt {2x})\right],$$
we will derive an ordinary differential equation, whose solution is $\Psi_{\tau}(x)$.  Precisely, it is given by

{\bf Lemma 2}.   {\it The kernel $\Psi_{\tau}(x)$ is a fundamental solution of the following fourth order differential equation with variable  coefficients}
$$x^2 {d^4 \Psi_{\tau} \over dx^4} + 5x {d^3 \Psi_{\tau} \over dx^3} + (4+\tau^2) {d^2 \Psi_{\tau} \over dx^2}
+ 16  \Psi_{\tau}  = 0. \eqno(1.26)$$

\begin{proof}  Recalling the Stirling asymptotic formula for the gamma-function \cite{erd},  Vol. I, we see that for each $\tau \in \mathbb{R}$ the integrand in (1.9) behaves as ($s= \gamma +it$)
$$\frac{ \Gamma(s/2)\Gamma((s+i\tau)/2) \Gamma((s-i\tau)/2)} {\Gamma((1-s)/2)} = 
e^{-\pi |t|/2} |t|^{2\gamma - 3/2},\quad  |t| \to \infty.$$
This argument  allows to differentiate repeatedly with respect to $x$ under the integral sign in (1.9).   Hence with the reduction formula for the gamma-function we obtain
$$x {d\over dx} x {d\over dx}\  \Psi_{\tau} =  - {1\over 16\pi i \sqrt\pi} \int_{\gamma-i\infty}^{\gamma+i\infty} \frac{  s^2 \Gamma(s/2)\Gamma((s+i\tau)/2) \Gamma((s-i\tau)/2)} {\Gamma((1-s)/2)} x^{-s} ds $$
$$= - \tau^2  \Psi_{\tau}(x)  - {1\over 4 \pi i \sqrt\pi} \int_{\gamma-i\infty}^{\gamma+i\infty} \frac{   \Gamma(s/2)\Gamma(1+ (s+i\tau)/2) \Gamma(1+ (s-i\tau)/2)} {\Gamma((1-s)/2)} x^{-s} ds $$
$$= -  \tau^2  \Psi_{\tau}(x)  - {1\over 4 \pi i \sqrt\pi} \int_{2+\gamma-i\infty}^{2+ \gamma+i\infty} \frac{   \Gamma(s/2- 1)\Gamma( (s+i\tau)/2) \Gamma( (s-i\tau)/2)} {\Gamma((3-s)/2)} x^{2-s} ds $$
$$=  -  \tau^2  \Psi_{\tau}(x)  - {1\over \pi i \sqrt\pi} \int_{2+\gamma-i\infty}^{2+ \gamma+i\infty} \frac{   \Gamma(s/2)\Gamma( (s+i\tau)/2) \Gamma( (s-i\tau)/2)} {(s-2)(1-s) \Gamma((1-s)/2)} x^{2-s} ds. $$
Differentiating again twice both sides of the latter equality and using (1.9),  we find
$${d^2\over dx^2} \left( x {d\over dx} \right) ^2 \Psi_{\tau} = -  \tau^2 {d^2 \Psi_{\tau} \over dx^2}  -  16 \Psi_{\tau}(x) . $$ 
Thus after fulfilling the differentiation in the left-hand side of the latter equality, we arrive at (1.26).  Lemma 2 is proved.

\end{proof}

Finally in this section we note that the obtained inequalities and integral representations of the kernel in (1.1)  will be used in the sequel to study the boundedness, compactness and invertibility of the index transforms (1.1), (1.2).

\section {Boundedness and compactness in  Lebesgue's spaces }

We begin,  introducing  the following weighted  $L_1$- space 
$$L^0\equiv   L_1\left(\mathbb{R}_+; \  K^2_{0}(\sqrt{2x})dx\right) := 
\left\{f:  \int_0^\infty K^2_{0}(\sqrt{2x})|f(x)|dx < \infty \right\}.\eqno(2.1)$$
In particular, as we will show below,  it contains  spaces  $L_{\nu, p}(\mathbb{R}_+)$ for some $\nu \in \mathbb{R},\   1\le  p \le \infty$ with the norms
$$||f||_{\nu, p}=\left(\int_0^\infty x^{\nu  p-1}|f(x)|^pdx\right)^{1/p}< \infty,\eqno(2.2)$$
$$||f||_{\nu, \infty}=\hbox{ess sup}_{x\ge 0}   |x^\nu f(x)|< \infty.$$
When $\nu={1\over p}$ we obtain the usual norm in $L_p$ denoted by   $|| \ ||_p$.  

{\bf Lemma 3.} {\it Let  $\  \nu <1, 1\le p \le \infty, \ q={p\over  p-1}$.   Then the embedding holds
$$L_{\nu, p}(\mathbb{R}_+)\subseteq L^0 \eqno(2.3)$$
and}
$$ ||f||_{L^0} \le 2^{-2(1/p+  \nu)} \left[ {\Gamma^{1/q}\left(2q(1-\nu)\right)\over
(2 q)^{2(1-\nu)} } \  B\left(1-\nu,\  1-\nu\right) \right]^2 \   ||f||_{\nu, p}, \ 1<p\le  \infty, \eqno(2.4)$$

$$ ||f||_{L^0} \le  \sup_{x\ge 0}  \left [K^2_{0}(\sqrt{2x} ) x^{1-\nu}\right] \   ||f||_{\nu,1}.\eqno(2.5)$$

\begin{proof}      In fact,  with the definition of the norm (2.1) and the H\" {o}lder inequality we obtain
$$||f||_{L^0} = \int_0^\infty K^2_{0}(\sqrt{2x}) |f(x)|dx \le \left(\int_0^\infty K_{0} ^{2q} (\sqrt{2x} )
x^{(1-\nu)q-1}dx\right)^{1/q}\  ||f||_{\nu, p},\   q={p\over  p-1}\eqno(2.6)$$
and the latter integral via asymptotic behavior of the Macdonald function (1.5), (1.6), (1.7) converges for 
$\nu  < 1$.   Hence integral (1.8) and  the generalized Minkowski  inequality yield
$$ \left(\int_0^\infty K_{0} ^{2q} (\sqrt{2x})
x^{(1-\nu)q-1}dx\right)^{1/q} = \left(\int_0^\infty x^{(1-\nu)q-1}
\left(\int_0^\infty e^{- \sqrt{2x} \  \cosh u} du\right)^{2q} dx\right)^{1/q}$$
$$\le \left(\int_0^\infty  \left(\int_0^\infty x^{(1-\nu)q-1} e^{- 2q \sqrt {2x} \  \cosh u}dx\right)
^{1/q} du \right)^2 $$$$ = 2^{2(\nu- 1/p)} (2 q)^{4(\nu-1)}\Gamma^{2/q}\left(2q(1-\nu)\right)\left( \int_0^\infty { du \over
\cosh^{2(1-\nu)} u}  \right)^2.$$
Calculating the integral with the hyperbolic function via relation (2.4.4.4) in \cite{prud},  Vol. I, we come up with the estimate (2.4).   For the case $p=1$ we end up  immediately  with (2.5), using (2.6), where the supremum is finite via the condition $\nu <1$.  Thus the embedding (2.3) is established and Lemma 3 is proved.

\end{proof}

{\bf Theorem 1.}   {\it The index transform  $(1.1)$  is well-defined as a  bounded operator from $L^0$ into the space $C_0(\mathbb{R})$ of bounded continuous functions vanishing at infinity.   Besides, the following composition representation holds
$$(F f)(\tau)=  \left( F_c (\mathcal{K}_0 f)( \cosh t)\right) (\tau),\eqno(2.7)$$
where the Fourier cosine transform $F_c$ is defined by $(1.14)$ and 
$$(\mathcal{K}_0  f)(x) =  -  \sqrt {{2\over \pi }}   \int_0^\infty    \ {\rm Re}  \left[ K_0\left( 4 e^{\pi i/4} \sqrt{ x t} \right) \right]\  f(t) dt\eqno(2.8)$$
is the operator of the Meijer type K- transform (cf. \cite{yal})}.

\begin{proof}  In fact, the inequality (1.25) implies  
$$||Ff||_{C_0(\mathbb{R})} = \sup_{\tau \in \mathbb{R}} |(F f) (\tau)|  \le {4\over \pi }  \int_0^\infty K^2_{0}(\sqrt {2x}) | f(x)| dx = {4\over \pi}  ||f||_{L^0} <  \infty ,$$
which means that the operator (1.1) is well-defined and the integral converges absolutely and uniformly with respect to $\tau \in \mathbb{R}$.  Thus $(F f) (\tau)$ is continuous.   On the other hand, recalling (1.15) and Corollary 3, we derive
$$|(F f) (\tau)|  \le {2\over \pi}  \int_0^\infty  \int_{0} ^\infty K_{0}\left(2\sqrt {2x\cosh t} \right)  |f(x)| dx dt \le 
{4\over \pi}   ||f||_{L^0}< \infty .$$
Hence in view of Fubini's theorem one  can invert the order of integration in the corresponding iterated integral and arrive at the composition (2.7).   Moreover,  the previous estimate says that $(\mathcal{K}_0 f)(\cosh t) \in L_1(\mathbb{R})$. Consequently, $(Ff)(\tau)$ vanishes at infinity owing to the Riemann-Lebesgue lemma. 
\end{proof} 

{\bf Corollary 4.}  {\it The operator $(1.1)$ \  $F : L_{\nu,p}(\mathbb{ R}_+) \to L_p(\mathbb{ R}), \ p \ge 2,  \   \nu < 1$ is bounded and }
$$||F f ||_{L_p(\mathbb{R})} \le  {\pi^{1/p- 1}\over 2}  \   q^{2(\nu-1)} \left[\Gamma\left(q(1-\nu)\right)\right]^{2/q} 
 B\left(1-\nu,\   1-\nu \right)  ||f||_{\nu,p}, \quad  q={p\over p-1}.\eqno(2.8)$$

\begin{proof}  Indeed,  taking the composition (2.7), we  appeal to the Hausdorff-Young inequality for  the Fourier cosine  transform (1.14) (cf. \cite{tit},  Theorem 74)
$$|| F_c f ||_{L_q(\mathbb{R})} \le (2\pi)^{1/q-1/2}  || f ||_{L_p(\mathbb{R})},\  1 < p \le 2,
\ q= {p\over p-1},$$
we find

$$||F f ||_{L_p(\mathbb{R})}  \le \sqrt 2 \  \pi^{1/p-1/2}   \left(\int_{0} ^\infty 
\left|(\mathcal{K}_0  f)(\cosh t)\right|^q dt \right)^{1/q}.\eqno(2.9)$$
Hence by the generalized Minkowski and H\"older inequalities  with relation  (2.16.2.2) from \cite{prud},  Vol. II  we obtain 
$$\sqrt 2 \  \pi^{1/p-1/2}   \left(\int_{0} ^\infty  \left|(\mathcal{K}_0  f)(\cosh t)\right|^q dt \right)^{1/q}
\le    2 \pi^{1/p- 1} \int_0^\infty   |f(x)|\left(\int_0^\infty   K^q_0\left( 2\sqrt{ 2x}\   \cosh( t/2) \right) dt\right)^{1/q}dx$$
$$\le 2^{1/q+ 1}  \pi^{1/p-1}    \int_0^\infty   \int_0^\infty  |f(x)|  \left(\int_0^\infty  e^{-2q \sqrt {2x} \cosh t \cosh u} dt\right)^{1/q} du dx$$
$$=   2^{1/q+ 1}  \pi^{1/p- 1}   \int_0^\infty   \int_0^\infty  |f(x)| \  K_0^{1/q} \left(2q \sqrt{2x} \cosh u\right) du dx$$
$$\le  2^{1/q+ 1}  \pi^{1/p- 1}     ||f||_{\nu,p}  \int_0^\infty   \left( \int_0^\infty   x^{(1-\nu)q-1} K_0 \left(2q\sqrt {2x} \ \cosh u\right) dx \right)^{1/q} du $$
$$=  2^{\nu }  \pi^{1/p-1}  q^{2(\nu-1)} \left[\Gamma\left(q(1-\nu)\right)\right]^{2/q} 
 ||f||_{\nu,p}  \int_0^\infty   \frac{du}{\cosh^{2(1-\nu)} u } $$
$$=  2^{-\nu}  \pi^{1/p-1}  q^{2(\nu-1)} \left[\Gamma\left(q(1-\nu)\right)\right]^{2/q} 
 B\left(1-\nu,\   1-\nu \right)  ||f||_{\nu,p}.$$
Consequently,  combining with (2.9),  we get (2.8).
\end{proof}

Now we investigate the compactness of the operator (1.1).

{\bf Theorem 2.} {\it The operator  $(1.1)$ \ $F  :  L_{\nu,p}(\mathbb{ R}_+) \to L_q(\mathbb{R}), \  1< p \le 2,   \nu < 1, \  q=p/(p-1)$ is compact}.

\begin{proof} The proof is based on approximation of the operator (1.1) by a sequence of compact operators of a finite rank with continuous kernels of compact support.  But to achieve this goal, it is sufficient to verify the following Hilbert-Schmidt-type condition
$$ \int_0^\infty\int_{-\infty}^\infty \left |\Psi_{\tau}(x)\right|^{q}
x^{(1-\nu)q-1}d\tau dx < \infty.\eqno(2.10)$$
In fact, similarly as above  we recall (1.15) and the generalized Minkowski inequality to deduce
$$\left( \int_0^\infty\int_{-\infty}^\infty \left |\Psi_{\tau}(x) \right|^{q} 
x^{(1-\nu)q-1}d\tau dx \right)^{1/q} $$$$\le  2^{1/2+ 1/p}  \  \pi^{1/q-1/2} \left( \int_0^\infty x^{(1-\nu)q-1} \left( \int_{0}^\infty 
K^p _{0}(2\sqrt{2x} \cosh  t ) dt \right)^{1/ (p-1)} dx   \right)^{1/q} $$
$$\le   2^{1/2+ 1/p}  \  \pi^{1/q-1/2}  \left( \int_0^\infty x^{(1-\nu)q-1} \left( \int_{0}^\infty   
\  K^{1/p} _{0}(2p \sqrt{2x}  \cosh t)\  dt \right)^q  dx   \right)^{1/q} $$
$$\le     2^{1/2+ 1/p}  \  \pi^{1/q-1/2} \int_0^\infty   \left( \int_0^\infty x^{(1-\nu)q-1} 
\  K^{1/(p-1)} _{0}(2p\sqrt {2x}  \cosh t)\  dx\right)^{1/q}  dt $$
$$=  2^{\nu -3/2}  \  \pi^{1/q-1/2} p^{2(\nu-1)} B \left( 1-\nu,\ 1-\nu\right) $$
$$\times  \left( \int_0^\infty x^{2(1-\nu)q-1}  \  K^{1/(p-1)} _{0}(x )\  dx\right)^{1/q}  \le   2^{\nu -3/2}  \  \pi^{1/q-1/2} p^{2(\nu-1)} B \left( 1-\nu,\ 1-\nu\right)$$
$$\times    \left( \int_0^\infty du \left(\int_0^\infty x^{2(1-\nu)q-1}  e^{-x\cosh u /(p-1)} dx\right)^{p-1} \right)^{1/p}  =   2^{- \nu +1/2- 2/p}  \  \pi^{1/q-1/2} q^{2(\nu-1)}   \Gamma^{1/q} (2q(1-\nu))$$
$$\times B \left( 1-\nu,\ 1-\nu\right)   B^{1/p} \left(p(1-\nu),\ p(1-\nu)\right) < \infty.$$

\end{proof}

Another representation of the transform (1.1) can be given via the Parseval equality for the Mellin transform (1.13) and the Mellin-Barnes integral representation (1.9). In fact, an immediate consequence of Theorem  87 in \cite{tit} and Stirling's asymptotic formula for the gamma-function is 

{\bf Theorem 3}. {\it Let $f \in L_{1-\nu, p}(\mathbb{R}_+),\ 1 < p \le 2.$  Then for all $\tau \in \mathbb{R}$}
$$(F f)(\tau) =   - {1\over 16\pi i \sqrt\pi} \int_{\nu -i\infty}^{\nu +i\infty} \frac{ \Gamma(s/2)\Gamma((s+i\tau)/2) \Gamma((s-i\tau)/2)} {\Gamma((1-s)/2)} f^*(1-s)  ds.\eqno(2.11)$$

Finally in this section we investigate the existence and boundedness of the  operator (1.2), which is  the adjoint of (1.1). In fact, it follows from  the general operator theory.  However, we will prove it directly, getting an explicit estimation of its norm.  Assuming $g(\tau) \in L_p(\mathbb{R}),\ 1 <p \le 2$ and recalling  the asymptotic formula (1.5) for the Macdonald function, we find that for each $x >0$ the function ${\rm Re}  \left[ K_0\left( 4 e^{\pi i/4} \sqrt{ x\cosh t} \right) \right]  \in  L_p(\mathbb{R}),\ 1 <p \le 2$.    Hence via the Parseval theorem for the Fourier transform (cf. \cite{tit}, Theorem 75) and equality (1.15),  operator (1.2) can be written as
$$(Gg)(x) = - \sqrt{{2\over \pi}} \int_{-\infty}^\infty  {\rm Re}  \left[ K_0\left( 4 e^{\pi i/4} \sqrt{ x\cosh t} \right) \right]  \left(\mathcal{F} g\right) (t)dt,\   x >0,\eqno(2.12)$$
where $\left(\mathcal{F} g\right) (t) \in L_q(\mathbb{R}),\  q= {p\over p-1}$ is the Fourier transform  of $g$ 
$$\left(\mathcal{F} g\right) (x) = {1\over \sqrt{2\pi} }\int_{-\infty}^\infty  g(t) e^{ixt} dt\eqno(2.13)$$
and the integral converges in the $L_p$-sense.

{\bf Theorem 4.}   {\it Let $g \in L_p(\mathbb{R}),\   1 < p \le 2$. Then operator  $(1.2)$  is well-defined  and for all $x >0$}
$$|(Gg) (x)|  \le    2^{1/(4p)- 1}   p^{-1/(2p)}  x^{-1/(4p)} \  B\left({1\over 4p},\    {1\over 4p} \right) || g ||_{L_p(\mathbb{R})} .\eqno(2.14)$$

\begin{proof}   Indeed, taking  (2.12), we recall the H\"{o}lder inequality,  the Hausdorff-Young inequality for the Fourier transform (2.13) \cite{tit}  and the generalized Minkowski inequality to obtain 
$$|(G g) (x)| \le \sqrt{{2\over \pi}}  \left(\int_{-\infty}^\infty K^p_0(2\sqrt{ 2x}\cosh( t/2) ) dt \right)^{1/p}  || \mathcal{F} g ||_{L_q(\mathbb{R})}   \le   2 \pi^{1/q- 1}  || g||_{L_p(\mathbb{R})} $$
$$\times  \int_{0}^\infty  \left( \int_{-\infty}^\infty  e^{- p\sqrt{2x} \   t^2  \  \cosh u} dt \right)^{1/p}  du =  \   2^{1-1/(4p)}   p^{-1/(2p)}  x^{-1/(4p)} \   || g ||_{L_p(\mathbb{R})} \int_{0}^\infty {du  \over \cosh^{1/(2p)} u }$$
$$=    2^{1/(4p)- 1}   p^{-1/(2p)}  x^{-1/(4p)} \  B\left({1\over 4p},\    {1\over 4p} \right) || g ||_{L_p(\mathbb{R})} , $$
which proves (2.14). 

\end{proof} 

{\bf Theorem  5.}  {\it The operator $G : L_p(\mathbb{ R}) \to L_{\nu,r}(\mathbb{ R}_+), \   1< p \le 2,  \  q= p/(p-1),    r \ge 1,  \    \nu >0 $ is bounded and }
$$||G g ||_{\nu,r} \le 2^{\nu -1 +1/r- 2/p } \pi^{1/q- 1}  { \Gamma^{1/r} (2\nu r)\Gamma^{2/p} \left( \nu p\right) \over r^{2\nu} \Gamma^{1/p} (2\nu p)} \    B \left( \nu,\   \nu\right)  \    || g ||_{L_p(\mathbb{R})},$$
where  $p, r$ have no dependence.

\begin{proof}  Again with  (2.12),  the generalized Minkowski,   H\"{o}lder inequalities and  the Hausdorff-Young inequality  for the Fourier transform  (2.13) we  find
$$||G g||_{\nu,r}  \le \sqrt {{2\over \pi}} \int_{-\infty}^\infty \left|  \left(\mathcal{F} g\right) (t)\right|  \left(\int_{0} ^\infty 
x^{\nu r -1} K^r_0(2\sqrt{2x} \cosh( t/2) ) dx \right)^{1/r}  dt$$$$\le 2^{1/p} \  \sqrt {{2\over \pi}}  || \mathcal{F} g||_{L_q(\mathbb{R})} \left( \int_{-\infty}^\infty  \left(\int_{0} ^\infty  x^{\nu r -1} K^r_0(2\sqrt{2x} \cosh t) dx \right)^{p/r}  dt\right)^{1/p}$$
$$\le   2^{1 -3\nu +1/r} \pi^{1/q- 1}     || g ||_{L_p(\mathbb{R})} \left( \int_{-\infty}^\infty  {dt\over \cosh^{2 \nu p} t} \right)^{1/p}  \left(\int_{0} ^\infty  x^{2\nu r -1} K^r_0(x) dx \right)^{1/r} $$
$$=     2^{1 - \nu +1/r- 2/p } \pi^{1/q- 1} \ B^{1/p} \left( \nu p,\   \nu p\right)  \left(\int_{0} ^\infty  x^{2\nu r -1} K^r_0(x) dx \right)^{1/r}  || g ||_{L_p(\mathbb{R})}  $$
$$\le   2^{1 - \nu +1/r- 2/p } \pi^{1/q- 1} \ B^{1/p} \left( \nu p,\   \nu p\right)  
\int_0^\infty  \left(\int_{0} ^\infty  x^{2\nu r -1} e^{-xr \cosh u} dx  \right)^{1/r} du  \   || g ||_{L_p(\mathbb{R})}  $$
$$=   2^{1 - \nu +1/r- 2/p } \pi^{1/q- 1}  r^{-2\nu} \Gamma^{1/r} (2\nu r)  B^{1/p} \left( \nu p,\   \nu p\right)  
\    || g ||_{L_p(\mathbb{R})} \int_0^\infty  {du  \over \cosh^{2\nu} u }$$
$$=   2^{\nu -1 +1/r- 2/p } \pi^{1/q- 1}  r^{-2\nu} \Gamma^{1/r} (2\nu r)  B^{1/p} \left( \nu p,\   \nu p\right)    B \left( \nu,\   \nu\right)  \    || g ||_{L_p(\mathbb{R})}. $$ 

\end{proof}

\section{Inversion theorems}   

The composition representation (2.11) and the properties of the Fourier and Mellin transforms are key ingredients to prove the inversion theorem for the index transform (1.1) $(Ff)(\tau)$.   Namely, we have

{\bf Theorem 6}.  {\it  Let $ f (t) \in L_{1-\nu, p}(\mathbb{R}_+) \cap L_1((1,\infty); t dt),\ 1 < p \le 2,  0< \nu < 1,\  q= p/(p-1), $    and let the Mellin transform $(1.11)$ satisfy the condition $f^*(1)=0$.  If, besides,  $\tau e^{\pi \tau}  F(\tau) \in L_1(\mathbb{R}_+)$, then for all $x> 0$ the following inversion formula holds}
$$f(x) =  - { 4\over  \pi  }  {d\over dx}  \int_{0} ^\infty  {\tau (Ff)(\tau) \over \cosh(\pi\tau/2)}   K_{i\tau}\left(2 \sqrt{2x} \right)\   {\rm Re }\ J_{i\tau} \left(2\sqrt{2x} \right) d\tau. \eqno(3.1)$$

\begin{proof}   In fact, since $ f \in L_{1-\nu, p}(\mathbb{R}_+)$ then by virtue of Theorem 86 in \cite{tit} its Mellin transform $f^*(s) \in L_q(1-\nu- i\infty, 1-\nu+i\infty)$. Hence it is not difficult to verify with the use of H\"{o}lder's inequality that integral (2.11) converges absolutely. Moreover, taking the Fourier cosine transform (1.14) of both sides of this equality, we change the order of integration in its  right-hand side by Fubini's theorem. Indeed, this is  because the absolute convergence of the corresponding iterated integral can be immediately verified,  employing inequalities  (1.18), (1.19).   Then,  recalling  (1.16), we come up with the equality 
$$\int_0^\infty (F f)(\tau) \cos(\tau u) d\tau =   - {1\over 8 i \sqrt\pi} \int_{\nu -i\infty}^{\nu +i\infty} \frac{ \Gamma(s/2)\Gamma(s) } {\Gamma((1-s)/2)} f^*(1-s)  \left( 2\cosh u \right)^{-s} ds.\eqno(3.2)$$
Further,  differentiating both sides of (3.2) with respect to $u$, which is definitely allowed under conditions of the theorem, we obtain
$$\int_0^\infty \tau (F f)(\tau) \sin (\tau u) d\tau =   - {\tanh  u \over 4 i \sqrt\pi} \int_{\nu -i\infty}^{\nu +i\infty} \frac{ \Gamma(1+ s/2)\Gamma( s) } {\Gamma((1-s)/2)} f^*(1-s)  \left( 2\cosh u \right)^{-s} ds.\eqno(3.3)$$
Meanwhile, after the substitution $p=  \cosh u$ in (3.4), we find  
$$- {2\over \sqrt \pi} \int_0^\infty \tau (F f)(\tau) {\sin \left(\tau \log\left(p +\sqrt{p^2-1}\right)\right)\over \sqrt{p^2-1}} d\tau = {1\over 2p\pi i} \int_{\nu -i\infty}^{\nu +i\infty} \frac{ \Gamma(1+ s/2)\Gamma( s) } {\Gamma((1-s)/2)} f^*(1-s)  \left( 2p\right)^{-s} ds.\eqno(3.4)$$
In the meantime, relations (2.16.2.2),  (2.16.6.1) in \cite{prud},  Vol. II  give the  representation
$$ {\sin\left(\tau \log\left(p + \sqrt{p^2-1} \right)\right)\over \sqrt{p^2-1} } =  {\sinh(\pi\tau) \over 4 \pi^2 i}  $$
$$\times   \int_{\nu - i\infty}^{\nu+i\infty}  \Gamma(s)  \  \Gamma\left({1-s+ i\tau\over 2}\right)  \Gamma\left({1-s- i\tau\over 2}\right) (2p)^{-s} ds,\ 0 < \nu < 1.\eqno(3.5)$$
Hence, substituting the latter expression in the left-hand side of (3.4) and changing the order of integration with the use of the  condition  $\tau e^{\pi \tau}  F(\tau) \in L_1(\mathbb{R}_+)$,  we derive after simple changes of variables
$${1\over 2\pi  \sqrt \pi}  \int_{1-\nu - i\infty}^{1-\nu+i\infty}  \Gamma(1-w) (2p)^{w}  \int_0^\infty \tau (F f)(\tau) \Gamma\left({w+ i\tau\over 2}\right)  \Gamma\left({w- i\tau\over 2}\right) d\tau \ dw $$
$$= \int_{-\nu -i\infty}^{-\nu +i\infty} \frac{ \Gamma(1- s/2)\Gamma( 1-s) } {\Gamma((1+s)/2)} f^*(1+s)  \left( 2p\right)^{s} {ds\over s}.\eqno(3.6)$$
On the other hand, under conditions of the theorem $f^*(1+s)$ is analytic and bounded in the strip $-\nu  < {\rm Re} s < 1$.  In fact, it follows from the estimate $(s= \mu +iy)$
$$|f^*(1+s)| \le \int_0^1 |f(t)| t^\mu dt + \int_1^\infty |f(t)| t^\mu dt \le  \left(\int_0^1 |f(t)|^p  t^{(1-\nu)p-1} dt\right)^{1/p} 
 \left(\int_0^1 t^{(\mu+\nu)q -1} dt\right)^{1/q}$$
 $$   + \int_1^\infty |f(t)| t^\mu dt \le (q(\mu+\nu ))^{-1/q} ||f||_{1-\nu,p} +  \int_1^\infty |f(t)| \  t dt  < \infty.$$
Therefore, since the integrand in the right-hand side of (3.6) tends to zero when $|{\rm Im} s| \to \infty$ in this strip and $f^*(1)=0$, one can shift the contour to the right via the Cauchy theorem, integrating over the vertical line 
$(1-\nu- i\infty, 1-\nu+i\infty)$.  Hence,
$${1\over 2\pi  \sqrt \pi}  \int_{1-\nu - i\infty}^{1-\nu+i\infty}  \Gamma(1-w) (2p)^{w}  \int_0^\infty \tau (F f)(\tau) \Gamma\left({w+ i\tau\over 2}\right)  \Gamma\left({w- i\tau\over 2}\right) d\tau \ dw $$
$$= \int_{1-\nu -i\infty}^{1-\nu +i\infty} \frac{ \Gamma(1- s/2)\Gamma( 1-s) } {\Gamma((1+s)/2)} f^*(1+s)  \left( 2p\right)^{s} {ds\over s}$$
and the uniqueness theorem for the inverse Mellin transform (1.12) of integrable functions \cite{tit} implies 
$${1\over 2\pi  \sqrt \pi}   \int_0^\infty \tau (F f)(\tau) \Gamma\left({s+ i\tau\over 2}\right)  \Gamma\left({s- i\tau\over 2}\right) d\tau =  \frac{ \Gamma(1- s/2) } {\Gamma((1+s)/2)} {f^*(1+s) \over s},\ s \in  (1-\nu- i\infty, 1-\nu+i\infty).$$
Thus,
$$ {f^*(1+s) \over s} =  {1\over 2\pi  \sqrt \pi}   \int_0^\infty \tau (F f)(\tau) \frac {\Gamma((1+s)/2)} { \Gamma(1- s/2) }
 \Gamma\left({s+ i\tau\over 2}\right)  \Gamma\left({s- i\tau\over 2}\right) d\tau,\eqno(3.7)$$
Further, employing the relation (8.4.23.11) in \cite{prud}, Vol. III, we find the value of the Mellin-Barnes integral 
$${1\over 2\pi i}  \int_{1-\nu - i\infty}^{1-\nu+i\infty} \frac {\Gamma((1+s)/2)} { \Gamma(1- s/2) }
 \Gamma\left({s+ i\tau\over 2}\right)  \Gamma\left({s- i\tau\over 2}\right)  s x^{-s} ds =- {8\sqrt \pi\over \cosh(\pi\tau/2)}$$
$$\times\   x {d\over dx} \left[ K_{i\tau}\left(2 \sqrt{2x} \right)\   {\rm Re }\ J_{i\tau} \left(2\sqrt{2x} \right) \right],\ x >0.\eqno(3.8)$$
Hence, taking the inverse Mellin transform in (3.7), changing the order of integration by Fubini's theorem, which is permitted by the imposed conditions,  we  arrive  at the inversion formula (3.1), completing the proof of Theorem 6.

\end{proof} 

Considering  the index transform (1.2),  it has 

{\bf Theorem 7}.  {\it  Let $g(z/i)$ be an even analytic function in the strip $D= \left\{ z \in \mathbb{C}: \ |{\rm Re} z | < \alpha < 1\right\} ,\  g(0)=g^\prime (0)=0, \  g(z/i)= o(1),\  |{\rm Im} z |\to \infty$ uniformly in $D$  and $g(x)  \in L_1(\mathbb{R})$.  Then,  if the index transform $(1.2)$ satisfies the condition $x^{-\gamma}  {d (Gg)/ dx} \in L_1(0,1), \  1/2 < \gamma < 1$, the following  inversion formula, which is written in terms of the Stieltjes integral holds valid for all $x \in \mathbb{R}$}  
$$  g(x) =  {4\over \pi}  \  x\sinh\left({\pi x\over 2}\right) \int_0^\infty  K_{ix}\left(2 \sqrt{2y} \right)\   {\rm Re }\ J_{ix} \left(2\sqrt{2y} \right) d (Gg)(y) . \eqno(3.9)$$

\begin{proof}   Indeed, substituting in (1.2) the expression of its  kernel  in terms of the Mellin - Barnes integral (1.9), we change the order of integration by Fubini's theorem via the condition $g(x)  \in L_1(\mathbb{R})$ and the estimate 
 $$\int_{-\infty}^\infty    |g(\tau)| \  \int_{\gamma -i\infty}^{\gamma +i\infty} \left| \frac{ \Gamma(s/2)\Gamma((s+i\tau)/2) \Gamma((s-i\tau)/2)} {\Gamma((1-s)/2)} ds \right| d\tau$$
$$\le B(\gamma/2,\gamma/2) \int_{-\infty}^\infty    |g(\tau)| d\tau \  \int_{\gamma -i\infty}^{\gamma +i\infty} \left| \frac{ \Gamma(s/2) \Gamma(s) } {\Gamma((1-s)/2)} ds \right| < \infty,\  \gamma  >0.$$
Hence
$$(Gf)(x)=  - {1\over 16 \pi i \sqrt \pi}  \int_{\gamma -i\infty}^{\gamma +i\infty}  \frac{ \Gamma(s/2)} {\Gamma((1-s)/2)} x^{-s} \int_{-\infty}^\infty    \Gamma\left({s+i\tau\over 2}\right) \Gamma\left({s-i\tau\over 2}\right) g(\tau) \  d\tau ds$$
and since $(Gf)(x) x^{\gamma-1} \in L_1(\mathbb{R}_+)$, which can be verified by moving the contour $(\gamma -i\infty,\  \gamma +i\infty)$ in the right-hand side of the previous equality to the right, the Mellin transform $(G g)^*(s)$ exists (see \cite{tit})  and can be represented reciprocally  in the form
$$(G g)^*(s)=  - {1\over 8\sqrt \pi} \  \frac{ \Gamma(s/2)} {\Gamma((1-s)/2)}  \int_{-\infty}^\infty    \Gamma\left({s+i\tau\over 2}\right) \Gamma\left({s-i\tau\over 2}\right) g(\tau) \  d\tau.\eqno(3.10)$$
Meanwhile,  the Stirling asymptotic formula for the gamma function yields  $(s=\gamma + iu)$
$$\frac {\Gamma((1-s)/2)} {   \Gamma(s/2)}  =O\left( |s|^{1/2- \gamma}\right),\  
|s| \to \infty.\eqno(3.11)$$
Moreover, from the definition of the index transform (1.2) and integral representation (1.9) of its kernel it is not difficult to verify that $(Gg)(x)$ is differentiable and $(Gg)^\prime (x) x^\gamma,\  \gamma >1/2$ vanishes at infinity. Hence, integrating by parts, we  write  the Mellin transform $(Gg)^*(s)$ as
$$(Gg)^*(s) = \int_0^\infty (Gg)(x) x^{s-1} dx =  - {1\over s} \int_0^\infty (Gg)^\prime (x) x^{s} dx,\ {\rm Re\  s} =\gamma > 0,\eqno(3.12)$$
we see that under condition $(Gg)^\prime (x) x^{-\gamma} \in L_1(0,1), \  1/2 < \gamma  < 1$ 
$$ 2^{s}\   \frac {\Gamma((1-s)/2)} { \Gamma(s/2)}  (G g)^*(s) =  - 2^{s-1}\   \frac {\Gamma((1-s)/2)} { \Gamma(1+ s/2)} 
 \int_0^\infty (Gg)^\prime (x) x^{s} dx \in L_1(\gamma-i\infty,\ \gamma +i\infty)$$ 
is analytic in the vertical strip $-\gamma < {\rm Re s} < \gamma$ and via (3.11) tends to zero when $|{\rm Im s}| \to \infty $ uniformly in the strip.  Appealing to  relation (2.16.2.2) in \cite{prud}, Vol. II,  the inverse Mellin transform (1.12) implies from (3.10)  
$$ {1\over i \sqrt \pi } \int_{\gamma -i\infty}^{ \gamma +i\infty}  2^{s}\   \frac {\Gamma((1-s)/2)} { \Gamma(s/2)}  (G g)^*(s) y^{-s} ds =  -  \int_{-\infty}^\infty    K_{i\tau}  (y) g(\tau) \  d\tau,\  y >0.\eqno(3.13)$$
Further,  writing $K_{i\tau} (y)$  in terms of the modified Bessel function of the third kind $I_{i\tau}(y)$ \cite{erd}, Vol. II, we get 
$$K_{i\tau}(y) = {\pi\over 2 i \sinh(\pi\tau)} \left[I_{-i\tau}(y)-  I_{i\tau}(y)\right].$$
Substituting it in the right-hand side of (3.13) and taking into account the evenness of $g(\tau)$, we find the equality
$$  {1\over i \sqrt \pi } \int_{\gamma -i\infty}^{ \gamma +i\infty}  2^{s}\   \frac {\Gamma((1-s)/2)} { \Gamma(s/2)}  (G g)^*(s) y^{-s} ds =    -\pi i  \int_{-i\infty}^{i\infty}     I_{z}  (y) g(z/i) \  {dz \over \sin(\pi z)}.\eqno(3.14)$$
 On the other hand, according to our assumption $g(z/i)$ is analytic in the vertical  strip $0\le  {\rm Re}  z < \alpha$, tending to zero when $|{\rm Im} | \to \infty $ uniformly in the strip. Thus appealing to the inequality for the modified Bessel   function of the third kind  (see \cite{yal}, p. 93)
 $$|I_z(y)| \le I_{  {\rm Re} z} (y) \  e^{\pi |{\rm Im} z|/2},\   0< {\rm Re} z < \alpha,$$
one can move the contour to the right in the right-hand side  and, in turn,  to the left in the left-hand side of the equality  (3.14),  taking into account the  analytic properties of the corresponding integrand.   Hence  after a simple substitution we obtain  
$$  {1\over i \sqrt \pi } \int_{\gamma -i\infty}^{ \gamma +i\infty}  2^{-s}\   \frac {\Gamma((1+s)/2)} { \Gamma(- s/2)}  (G g)^*(-s) y^{s} ds =    -\pi i  \int_{\alpha -i\infty}^{\alpha+ i\infty}     I_{z}  (y) g(z/i) \  {dz \over \sin(\pi z)}.\eqno(3.15)$$
Multiplying both sides of (3.15) by $K_{ix} (y) y^{-1}$ and integrating over $(0,\infty)$, we interchange  the order of integration in both sides by the Fubini theorem. Then  employing relation (2.16.2.2) in \cite{prud}, Vol. II and the value of the integral (see relation (2.16.28.3) in \cite{prud}, Vol. II
$$\int_0^\infty K_{ix}(y) I_z(y) {dy\over y} = {1\over x^2 + z^2},$$ 
we come up with the equality 
$$  {1\over 4 i \sqrt \pi } \int_{\gamma -i\infty}^{\gamma +i\infty}   \frac {\Gamma((1+s)/2)} { \Gamma(-s/2)}   \Gamma\left({s+ix\over 2}\right) \Gamma\left({s-ix\over 2}\right) (G g)^*(-s)  ds =    -\pi i  \int_{\alpha -i\infty}^{\alpha+ i\infty}     {  g(z/i) \  dz \over (x^2+ z^2) \sin(\pi z)}.$$
Now, taking in mind (3.8), (3.12) and the Parseval identity  (1.13), the previous equality becomes 
$$  {1\over \cosh(\pi x/2)} \int_0^\infty  K_{ix}\left(2 \sqrt{2y} \right)\   {\rm Re }\ J_{ix} \left(2\sqrt{2y} \right) d (Gg)(y)  =   {i \over 2}  \int_{\alpha -i\infty}^{\alpha+ i\infty}     {  g(z/i) \  dz \over (x^2+ z^2) \sin(\pi z)}.\eqno(3.16)$$
Meanwhile, with the evenness of $g$ we write the right-hand side of (3.16) in the form
 $$ {i \over 2}  \int_{\alpha -i\infty}^{\alpha+ i\infty}     {  g(z/i) \  dz \over (x^2+ z^2) \sin(\pi z)}
 ={ i \over 4}  \left( \int_{-\alpha +i\infty}^{- \alpha- i\infty}   +   \int_{\alpha -i\infty}^{ \alpha+  i\infty}   \right)  {  g(z/i) \  dz \over (z-ix) \  z \sin(\pi z)} =  {\pi g(x) \over  2 x\sinh(\pi x)},$$
 appealing   to the Cauchy formula, because $ g(z/i) / ( z \sin(\pi z))$ is analytic in the strip $|{\rm Re} z| < \alpha < 1$ and tends to zero when $|{\rm Im} z| \to \infty$ uniformly in the strip. Thus, combining with (3.16),  we arrived at the inversion formula (3.9) and complete to proof of the theorem. 
 
 \end{proof}

\section{Initial   value problem}

The index transform (1.2) can be successfully applied to solve an initial value  problem  for the following fourth order partial differential  equation, involving the Laplacian
$$\left[ \left(x {\partial \over \partial x}  + y {\partial \over \partial y} \right)^2  + 
4 \left(x {\partial \over \partial x}  + y {\partial \over \partial y} \right)  +  2 \right]  \Delta u +  16  u  =0, \eqno(4.1)$$ 
where $\Delta = {\partial^2 \over \partial x^2} +  {\partial^2 \over \partial y ^2}$ is the Laplacian in $\mathbb{R}^2$.   Indeed, writing equation (4.1) in polar coordinates $(r,\theta)$, we end up with the equation 
$$ { \partial^2  \over \partial r^2}  \left(  \left(  r {\partial  \over \partial r} \right)^2 u  +    {\partial^2 u \over \partial \theta^2}\right) + 16 u = 0.\eqno(4.2)$$

{\bf Lemma 4.} {\it  Let $g(\tau)  \in L_1\left(\mathbb{R}; e^{ \beta |\tau|} d\tau\right),\  \beta \in (0, 2\pi)$. Then  the function
$$u(r,\theta)=  \int_{-\infty}^\infty    K_{i\tau}(2\sqrt {2r}) \ {\rm Im} \left[ J_{i\tau} (2\sqrt {2r}) \right] \frac{ e^{\theta\tau} g(\tau)} {\sinh(\pi\tau /2)} \  d\tau\eqno(4.3)$$
 satisfies   the partial  differential  equation $(4.2)$ on the wedge  $(r,\theta): r   >0, \  0\le \theta <  \beta$, vanishing at infinity.}

\begin{proof} In fact, the proof  follows from the direct substitution (4.3) into (4.2) and the use of Lemma 2.  The necessary  differentiation  with respect to $r$ and $\theta$ under the integral sign is allowed via the absolute and uniform convergence, which can be justified  with the use of (1.9) and the integrability condition $g \in L_1\left(\mathbb{R}; e^{ \beta |\tau|} d\tau\right),\  \beta \in (0, 2\pi)$ of the lemma.
\end{proof}

Finally, as a direct consequence of Theorem 7, we will formulate the initial value problem for equation (4.2) and give its solution.

{\bf Theorem 8.} {\it Let 
$$ g(x) =  {4\over \pi}  \  x\sinh\left({\pi x\over 2}\right) \int_0^\infty  K_{ix}\left(2 \sqrt{2y} \right)\   {\rm Re }\ J_{ix} \left(2\sqrt{2y} \right) d G (y)$$
and $G(y)$  satisfy conditions of Theorem $7$. Then  $u (r,\theta),\   r >0,  \  0\le \theta < \beta$ by formula $(4.3)$  will be a solution  of the initial value problem for the partial differential  equation $(4.2)$ subject to the initial condition}
$$u(r,0) = G (r).$$

\bigskip
\centerline{{\bf Acknowledgments}}
\bigskip

The work was partially supported by CMUP (UID/MAT/00144/2013), which is funded by FCT(Portugal) with national (MEC) and European structural funds through the programs FEDER, underthe partnership agreement PT2020.

\bibliographystyle{amsplain}

\begin{thebibliography}{10}

\bibitem{yak}   Yakubovich S.  Index transforms.   Singapore:  World Scientific Publishing Company; 1996.

\bibitem{erd}    Erd\'elyi A,  Magnus W,   Oberhettinger  F,   Tricomi FG.  Higher transcendental functions. Vols. I,  II. New  York: McGraw-Hill;  1953.

\bibitem{prud}  Prudnikov AP,  Brychkov  YuA,  Marichev OI. Integrals and series:  Vol. I: Elementary functions. New York:  Gordon and Breach;   1986;   Vol. II:  Special functions. New York: Gordon and Breach;  1986;   Vol. III:  More special functions. New York:   Gordon and Breach; 1990.


\bibitem{yal}   Yakubovich S,  Luchko Yu.  The hypergeometric approach to integral transforms and convolutions, Mathematics and its applications.  Vol. 287.  Dordrecht:  Kluwer Academic Publishers Group; 1994.

\bibitem {tit}  Titchmarsh EC.   An introduction to the theory of Fourier integrals.   New York:  Chelsea; 1986.





\end{thebibliography}

\end{document}